\RequirePackage[l2tabu, orthodox]{nag}
\documentclass[fontsize=14pt,headings=small,open=any,twoside=false,draft=false]{scrartcl}

\usepackage{cmap}
\usepackage{pdfsync}

\usepackage[T1,T2A]{fontenc}
\usepackage[cp1251]{inputenc}
\usepackage[french,english]{babel}

\usepackage[a4paper,pdftex,dvips,nohead,verbose,includehead,
            headheight=17pt,headsep=7mm,top=20mm,
            text={170mm,247mm},bottom=20mm,%%,
]{geometry}

\usepackage{titlesec}
\titleformat{\section}[hang]{\normalfont\large}{\bfseries\thesection.}{.5em}{\bfseries}
\titlespacing{\section}{\parindent}{12pt plus 6pt}{3pt}%{\wordsep}
\titleformat{\subsection}[runin]{\normalfont}{\bfseries\thesubsection.}{.5em}{\bfseries}[.\quad]
\titleformat{\subsubsection}[runin]{\normalfont}{\bfseries\thesubsection.}{.5em}{\bfseries}[.\quad]
\makeatletter
\renewcommand\l@section{\@dottedtocline{1}{0ex}{2em}}
\makeatother

\usepackage{fancyhdr}
\usepackage[format=plain,labelfont=it,font=normalsize,skip=6pt]{caption}
\DeclareCaptionLabelSeparator{dot}{.~} % разделитель между заголовком и телом подписи
\usepackage{amsmath,amsfonts,amssymb,latexsym,amsthm,multicol}
\usepackage{float}
\usepackage{indentfirst} % абзацный отступ в первом параграфе
\usepackage{graphicx}
\usepackage{wrapfig}
\usepackage{cite} % правильная организация ссылок в команде \cite
\usepackage{ncccomma}
\usepackage{soul}

\usepackage[shortlabels]{enumitem}
\usepackage[section,nottoc,notlot,notlof]{tocbibind}
\usepackage{mathrsfs}
\usepackage{multicol}
\usepackage{url}
\setlist{noitemsep}
\newcommand{\be}{\begin{eqnarray}}
\newcommand{\ee}{\end{eqnarray}}
\newcommand{\bew}{\begin{eqnarray*}}
\newcommand{\eew}{\end{eqnarray*}}

\newcommand{\Arn}{V.\,I.\,Arnold}
\newcommand{\Kozl}{V.\,V.\,Kozlov}

\begin{document}
%%%%%%%%%%%%%%%%%%%
% Пример титульной страницы
%%%%%%%%%%%%%%%%%%%

%\thispagestyle{empty}
\begin{center}
\bigskip
\textbf{\Large
\uppercase{On geometric methods\\ in works by\\[5pt] V.\,I.\,Arnold and V.\,V.\,Kozlov}\footnote{This work is supported by Russian Fund of Basic Research, Project No~11--01--00023-a.}
}

\bigskip
\textit{\large A.\,D.\,Bruno}

\bigskip
\textit{Keldysh Institute of Applied Mathematics, Moscow, Russia}
\end{center}

\begin{abstract}
We give a survey of geometric methods used in papers and books by \Arn\ and by \Kozl. They are methods of different normal forms, of different polyhedra, of small denominators and of asymptotic expansions.
\end{abstract}


\noindent

\section*{Introduction}
\addcontentsline{toc}{section}{Introduction}
In paper~\cite{cite0} there was given a short description of main achievements of \Arn. Below in \ref{sect1}--\ref{sect3} we give some additions to several Sections of this paper. In Sections~\ref{sect4}, \ref{sect5}, \ref{sect7}, \ref{sect8} we discuss two kinds  of normal forms in papers by \Arn\ and by \Kozl.

Logarithmic branching of solutions to Painlev\'e equations are discussed in~\ref{sect6}.

\section[On small divisors]{On the last paragraph of page~381~\cite{cite0} devoted to small divisors} \label{sect1}
Arnold's Theorem on stability of the stationary point in the Hamiltonian system with two degrees of freedom in~\cite{cite1} had wrong formulation (see~\cite[\S~12, Section IVd]{cite2}). Then \Arn~\cite{cite3} added one more condition in his Theorem, but its proof was wrong because it used the wrong statement (see~\cite{cite4, cite5}). All mathematical world was agreed with my critics except  \Arn. On the other hand, in the first proof of the same Theorem by J.\,Moser~\cite{cite6} there was a similar mistake (see~\cite[\S~12, Section~IVe]{cite2}). But in~\cite{cite7} J.\,Moser corrected his proof after my critics, published in~\cite[\S~12, Section~IVe]{cite2}.

Concerning the KAM theory. In 1974 I developed its generalization via normal forms~\cite[Part II]{cite8, cite9, cite10}. But up-to-day almost nobody understands my generalization.

\section[On higher-dimensional analogue of continued fraction]{On the last paragraph of page~384~\cite{cite0} concerning higher-dimensional analogue of continued fraction} \label{sect2}
The paper~\cite{cite11} \foreignlanguage{french}{``Poly\`edre d'Arnol'd et \ldots''} by G.~Lachaud (1993) was presented to C.R. Acad. Sci. Paris by \Arn. When I saw the article I published the paper~\cite{cite12} ``Klein polyhedrals \ldots'' (1994), because so-called ``Arnold polyhedra'' were proposed by F.\,Klein hundred years early. Moreover, they were introduced by B.\,F.\,Scubenko in 1988 as well. In 1993--2003 me and V.\,I.\,Parusnikov studied Klein polyhedra from algorithmical view point and found that they do not give a basis for algorithm generalizing the continued fraction. So I proposed another approach and another sole polyhedron, which give a basis for the generalization in 3 and any dimension (see~\cite{cite13, cite14, cite15, cite16, cite17}). Now there are a lot of publications on the Klein polyhedra and their authors following after \Arn\ wrongly think that the publications are on generalization of the continued fraction.

\section[On Newton Polygon]{On the last two paragraphs of page~395~\cite{cite0} devoted to Newton Polygon}\label{sect3}
In that text the term ``Newton polygon'' must be replaced by ``Newton polyhedron''. In contemporary terms I.\,Newton introduced \textit{support} and one \textit{extreme edge} of the \textit{Newton open polygon} for one polynomial of two variables. The full Newton open polygon was proposed by V.\,Puiseux (1850) and by C.~Briot and T.~Bouquet (1856) for one ordinary differential equation of the first order. Firstly a polyhedron as the convex hull of the support was introduced in my paper~\cite{cite18} (1962) for an autonomous system of $n$ ODEs. During 1960--1969 V.~I.~Arnold wrote 3 reviews on my works devoted to polygons and polyhedrons for ODEs with sharp critics ``of the geometry of power exponents'' (see my book~\cite[Ch.~8, Section~6]{cite19}). Later (1974) he introduced the name ``Newton polyhedron'', made the view that it is his invention and never gave references on my work. Now I have developed ``Universal Nonlinear Analysis'' which allows to compute asymptotic expansions of solutions to equations of any kind (algebraic, ordinary differential and partial differential)~\cite{cite20}.

\section{On non-Hamiltonian normal form}\label{sect4}
%There were other cases where Arnold was a \textbf{thief}. For example, 
In my candidate thesis ``Normal form of differential equations''~\cite{cite21} (1966) I introduced normal forms in the form of power series. It was a new class of them. Known before normal forms (NF) were either linear (Poincare, 1879)~\cite{cite22} or polynomial (Dulac, 1912)~\cite{cite23}. An official opponent was A.N. Kolmogorov. He estimated very high that new class of NF. Arnold put my NF into his book~\cite[\S~23]{cite24} without reference on my publication and named it as ``Poincare-Dulac normal form''. So, readers of his book attributed my NF to Arnold. I saw several articles where my NF were named as Arnold's.
%Indeed, V.~I.~Arnold made a lot of mistakes, but I do not know a case when he recognize his mistake.

\section{On canonical normalizing transformation}\label{sect5}
In~\cite[Ch.~7, \S~3, Subsection~3.1]{cite25} a proof of Theorem~7 is based on construction of a generating function $F=\left<P,q\right>+S_l(P,q)$ in mixed coordinates $P, q$. Transformation from old coordinates $P,Q$ to new coordinates $p, q$ is given by the formulae
\begin{equation}\label{eq1star}
p=\frac{\partial F}{\partial q},\quad Q=\frac{\partial F}{\partial P}.\tag{1*}
\end{equation}
Here $S_l(P,q)$ is a homogeneous polynomial in $P$ and $q$ of order $l$. According to~\eqref{eq1star}, the transformation from coordinates $P, Q$ to coordinates $p, q$ is given by infinite series, which are results of the resolution of implicit equations~\eqref{eq1star}. Thus, the next to the last sentence on page~272 (in Russian edition) ``The normalizing transformation is constructed in the form of a polynomial of order $L-1$ in phase variables'' is wrong. Indeed that property has the normalizing transformation computed by the Zhuravlev-Petrov method~\cite{cite26}.

\section{On branching of solutions of Painlev\'e equations}\label{sect6}
In~\cite[Ch.~I, \S~4, example~1.4.6]{cite27} the Painlev\'e equations are successive considered. In particularly, there was find the expansion
\begin{equation}\label{eq2star}
x(\tau)=\tau^{-1}\sum_{k=0}^\infty x_k\tau^k\tag{2*}
\end{equation}
of a solution to the fifth Painlev\'e equation. The series~\eqref{eq2star} is considered near the point $\tau=0$. After the substitution $\tau=\log t$, we obtain the series
\begin{equation}\label{eq3star}
x(t)=\log^{-1}t\sum_{k=0}^\infty x_k\log^kt,\tag{3*}
\end{equation}
which has a sense near the point $t=1$, where $\log t=0$. However, from the last expansion~\eqref{eq3star} authors concluded that $t=0$ is the point of the logarithmic branching the solution $x(t)$. It is wrong, because the expansion~\eqref{eq3star} does not work for $t=0$ as $\log0=\infty$ and the expansion~\eqref{eq3star} diverges.
That mistake is in the first edition of the book~\cite{cite27} (1996) and was pointed out in the paper~\cite{cite28} (2004), but it was not corrected in the second ``corrected'' edition of the book~\cite{cite27}.

A similar mistake is there in consideration of the sixth Painlev\'e equation. There for a solution to the sixth Painlev\'e equation, it was obtained the expansion~\eqref{eq2star}. After the substitution $\tau=\log(t(t-1))$, it takes the form
\[
x(t)=\log^{-1}(t(t-1))\sum_{k=0}^\infty x_k\log^k(t(t-1)).
\]
As the expansion~\eqref{eq2star} has a sense near the point $\tau=0$, the last expansion has a sense near points $t=(1\pm\sqrt5)/2$, because in them $t(t-1)=1$ and $\tau=0$. Thus, the conclusion in the book, that points $t=0$ and $t=1$ are the logarithmic branching points of the solution, is non correct. The mistake was point out in the paper~\cite{cite29} (2004), but it was repeated in the second edition of the book~\cite{cite27}. Indeed solutions of Painlev\'e equations have logarithmic branching, see~\cite{cite30, cite31}.

\section{On integrability of the Euler-Poisson equations}\label{sect7}
In the paper~\cite{cite32} Theorem~1 on nonexistence of an additional analytic integral was applied in \S~3 to the problem of motion of a rigid body around a fixed point. The problem was reduced to a Hamiltonian system with two degrees of freedom and with two parameters $x, y$. The system has a stationary point for all values of parameters. Condition on existence of the resonance $3:1$ was written as equation~(6) on parameters $x, y$. Then the second order form of the Hamiltonian function was reduced to the simplest form by a linear canonical transformation
\begin{equation}\label{eq4star}
(x_1,x_2,y_1,y_2)\to(q_1,q_2,p_1,p_2).\tag{4*}
\end{equation}

Condition of vanishing the resonant term of the fourth order in the obtained Hamiltonian function was written as equation~(7) on $x, y$. System of equations~(6) and~(7) was considered for
\[
x>0\text{ and } y>\frac x{x+1},
\]
where the system has two roots
\begin{equation}\label{eq5star}
x=\frac43,\; y=1\text{ and }x=7,\; y=2.\tag{5*}
\end{equation}
They correspond to two integrable cases $y=1$ and $y=2$ of the initial problem. It was mentioned in Theorem~3. But in the whole real plane $(x,y)$ the system of equations~(6) and~(7) has roots~\eqref{eq5star} and three additional roots
\begin{align}
x&=-\frac{16}3,& y&=1;&&&&& x&=-\frac{17}9,& y&=2;\label{eq6star}\tag{6*}\\
&&&&x&=0,&y&=9.&&&&\label{eq7star}\tag{7*}
\end{align}
Roots~\eqref{eq6star} belong to integrable cases $y=1$ and $y=2$. But the root~\eqref{eq7star} is out of them. Indeed the transformation~\eqref{eq4star} is not defined for $x=0$. If to make an additional analysis for $x=0$, then for resonance $3:1$ one obtains two points:~\eqref{eq7star} and
\begin{equation}\label{eq8star}
x=0,\quad y=\frac19.\tag{8*}
\end{equation}
In both these points, the resonant term of the fourth order part of the Hamiltonian function vanishes. But points~\eqref{eq7star} and~\eqref{eq8star} are out of the integrable cases $y=1$ and $y=2$; they contradict to statement of Theorem~3~\cite{cite32}. The paper~\cite{cite32} was repeated in the book~\cite[Ch.~VI, \S~3, Section~3]{cite33}. A non-Hamiltonian study of the problem see in the paper~\cite[Section~5]{cite34}.

\section{On normal forms of families of linear Hamiltonian systems}\label{sect8}
Real normal forms of families of linear Hamiltonian systems were given in~\cite[\S~2]{cite35}, where formula~(16) wrongly indicated the normal form corresponding to the elementary divisor $\lambda^{2l}$: the third sum in the formula~(16) has to be omitted. The indicated mistake was reproduced in the first three editions of the book~\cite[Appendix~6]{cite36} by Arnold. Discussions of that see in the paper~\cite{cite37} and in the book~\cite[Ch.~I, Section~6, Notes to Section~1.3]{cite38}.

\section{Conclusions}\label{sect9}
\ref{sect1}--\ref{sect3} were sent to Notices of the AMS for publication as a Letter to the Editor. But Editor S.\,G.\,Krantz rejected it. I consider that as one more case of the scientific censorship in the AMS.

%\newpage

%\noindent
%Alexander Bruno\\
%\null Professor and Head\\
%\null Department of Singular Problems\\
%\null Keldysh Institute of Applied Mathematics,\\
%\null Department of Theory of Dynamical Systems\\
%\null Lomonosov Moscow State University\\
%\null\url{http://en.wikipedia.org/wiki/Alexander_Dmitrievich_Bruno}\\
\end{document}